# DETERMINING THE DIMENSION OF ITERATIVE HESSIAN TRANSFORMATION

By R. Dennis Cook[1] and Bing Li[2]

*University of Minnesota and Pennsylvania State University*

The central mean subspace (CMS) and iterative Hessian transformation (IHT) have been introduced recently for dimension reduction when the conditional mean is of interest. Suppose that $X$ is a vector-valued predictor and $Y$ is a scalar response. The basic problem is to find a lower-dimensional predictor $\eta^T X$ such that $E(Y|X) = E(Y|\eta^T X)$. The CMS defines the inferential object for this problem and IHT provides an estimating procedure. Compared with other methods, IHT requires fewer assumptions and has been shown to perform well when the additional assumptions required by those methods fail. In this paper we give an asymptotic analysis of IHT and provide stepwise asymptotic hypothesis tests to determine the dimension of the CMS, as estimated by IHT. Here, the original IHT method has been modified to be invariant under location and scale transformations. To provide empirical support for our asymptotic results, we will present a series of simulation studies. These agree well with the theory. The method is applied to analyze an ozone data set.

**1. Introduction.** The basic problem of dimension reduction for regression [Li (1991, 1992), Cook and Weisberg (1991) and Cook (1998a)] is to find a lower-dimensional predictor that carries all the information relevant to the regression. Suppose that $X$ is a $p$-dimensional predictor and $Y$ is a scalar response. If there is a $p$ by $q$, $q \leq p$, matrix $\eta$ such that the $q$ linear combinations $\eta^T X$ fully describe the conditional distribution of $Y$ given $X$, then the subspace spanned by the columns of $\eta$ is called a dimension reduction subspace. In symbols, if

$$Y \perp\!\!\!\perp X | \eta^T X,$$

Received November 2002; revised January 2004.
[1]Supported in part by NSF Grant DMS-01-03983.
[2]Supported in part by NSF Grant DMS-02-04662.
*AMS 2000 subject classifications.* Primary 62G08; secondary 62G09, 62H05.
*Key words and phrases.* Dimension reduction, conditional mean, asymptotic test, order determination, eigenvalues.

then the column space of $\eta$ is a dimension reduction subspace. Here $\perp\!\!\!\perp$ stands for independence, so the statement is that $Y$ is independent of $X$ given $\eta^T X$. Any subspace that contains a dimension reduction subspace is itself a dimension reduction subspace. Under mild conditions the intersection of all dimension reduction subspaces is itself a dimension reduction subspace and then is called the central subspace (CS), and written as $\mathcal{S}_{Y|X}$. If the CS is known, then $X$ can be replaced with $P_{\mathcal{S}_{Y|X}}X$ without loss of information on the conditional distribution of $Y|X$, where $P_{(\cdot)}$ indicates a projection in the usual inner product onto the indicated subspace.

If, as in many regression analyses, the conditional mean $E(Y|X)$ is of particular interest, then it is possible, and beneficial, to carry out dimension reduction for that purpose. Cook and Li (2002) formulated dimension reduction in this context as follows. If there is a $p$ by $q$ matrix $\eta$ such that

$$E(Y|X) = E(Y|\eta^T X),$$

then the column space of $\eta$ is a dimension reduction subspace for the conditional mean, and is called a *mean subspace*. Under mild conditions, the intersection of all such subspaces is again a mean subspace and then is called the central mean subspace (CMS), and written as $\mathcal{S}_{E(Y|X)}$. See Li, Cook and Chiaromonte (2003) and Yin and Cook (2002) for other developments related to the CMS.

Several benefits accrue from studying the conditional mean $E(Y|X)$ rather than all of $Y|X$:

1. Because $\mathcal{S}_{E(Y|X)} \subseteq \mathcal{S}_{Y|X}$, it may be possible to achieve further reduction of dimension.
2. As in classical estimation, focusing on a smaller inferential object could lead to increased accuracy. Here $\mathcal{S}_{E(Y|X)}$ acts as the "parameter of interest" and all aspects of the conditional distribution of $Y|X$ not described by the conditional mean act as the nuisance parameter.
3. Study of $\mathcal{S}_{E(Y|X)}$ leads to a categorization of several existing methods, such as ordinary least square estimates (OLS) [Li and Duan (1989)], sliced inverse regression (SIR) [Li (1991)], principal Hessian directions (PHD) [Li (1992)] and the sliced average variance estimator (SAVE) [Cook and Weisberg (1991)]. It thus provides further insight into the dimension reduction problem.

As demonstrated by Cook and Li (2002), these four methods estimate either the CS or the CMS under the first or both of the following two conditions:

(A) *Linearity condition*: $E(X|P_\mathcal{S} X)$ is a linear function of $X$,
(B) *Constant covariance condition*: $\text{Var}(X|P_\mathcal{S} X)$ is a nonrandom matrix,

where the subspace $\mathcal{S}$ is either the CS $\mathcal{S}_{Y|X}$ or the CMS $\mathcal{S}_{E(Y|X)}$, depending on the method. In particular, using $\mathcal{S} = \mathcal{S}_{E(Y|X)}$, OLS and PHD estimate

vectors in the CMS, with OLS requiring condition (A) and PHD requiring both conditions. Using $\mathcal{S} = \mathcal{S}_{Y|X}$, SIR and SAVE estimate vectors in the CS, with SIR requiring condition (A) and SAVE requiring both conditions.

Condition (A) holds for all subspaces of $\mathbb{R}^p$ if the predictor $X$ has an elliptical distribution [Eaton (1986)], and it holds approximately if $\dim(\mathcal{S}) \ll p$ [Hall and Li (1993)]. Condition (B) is more stringent but will be satisfied if $X$ has a multivariate normal distribution. It is noteworthy that both conditions apply to the marginal distribution of the predictors and not to the conditional distribution of $Y|X$ as is common in regression modeling. Consequently, we are free to use experimental design, predictor transformations or re-weighting [Cook and Nachtsheim (1994)] to induce the conditions as necessary without suffering complications when inferring about $E(Y|X)$ or $Y|X$.

In some practical problems, such as in the recumbent cow data [Clark et al. (1987); see also Cook and Li (2002)], there is significant heteroscedasticity among the predictors. In such cases condition (B) fails and application of PHD and SAVE becomes problematic. Also, OLS provides at most one vector, and so does SIR if the response $Y$ is binary as it is in the recumbent cow data. Hence if the dimension of the CMS or the CS is 2 or more, then OLS and SIR will necessarily miss part of the CMS and CS.

It is in this context that Cook and Li (2002) introduced the method of IHT, presenting two versions of IHT that both estimate vectors in the CMS. Assuming that $X$ is standardized to have mean zero and covariance matrix $I_p$ (the identity matrix of dimension $p$), one version uses the response-based (y-based) Hessian matrix $E((Y - E(Y))XX^T)$ and the other uses the residual-based (r-based) Hessian matrix

$$H = E((Y - E(Y) - E(YX^T)X)XX^T).$$

Both versions require only condition (A) but, like PHD, can estimate multiple vectors in the CMS. Following the findings of Cook (1998b) on the general superiority of r-based PHD over y-based PHD, we use the r-based Hessian matrix $H$ in the rest of this article.

Cook and Li (2002) demonstrated the following fundamental relation, which is the basis for IHT. Under condition (A) alone, the CMS is an invariant subspace of the linear transformation $H$, that is,

$$H\mathcal{S}_{E(Y|X)} \subseteq \mathcal{S}_{E(Y|X)}.$$

It follows that if we know any nonzero vector in the CMS, then we can transform it iteratively by the Hessian matrix to bring out other vectors in the CMS. An obvious initial vector is the OLS vector, which we know belongs to the CMS under condition (A) alone. Thus if we use $\beta$ to denote the OLS vector, which, assuming $X$ to be standardized, has the form $\beta = E(YX)$, then the vectors

$$\beta, H\beta, H^2\beta, \ldots$$

are all in the CMS. At a certain point, say at $H^{k-1}\beta$ (with $k \leq p$), one more iteration ceases to bring out a linearly independent vector, and all subsequent vectors in the sequence must also be linear combinations of the first $k$ vectors.

In brief, under the linearity condition (A) the *IHT subspace*

$$\mathcal{S}_{\text{IHT}} \equiv \text{Span}\{\beta, H\beta, H^2\beta, \ldots, H^{p-1}\beta\} = \text{Span}\{\beta, H\beta, H^2\beta, \ldots, H^{k-1}\beta\}$$

is contained in the CMS, $\mathcal{S}_{\text{IHT}} \subseteq \mathcal{S}_{E(Y|X)}$ [Cook and Li (2002)]. The task of this paper is to estimate the dimension $k$ of $\mathcal{S}_{\text{IHT}}$, the number of linearly independent vectors generated by the iterative transformations $H^j\beta$, $j = 0, 1, \ldots, p-1$, and thereby obtain an estimator of $\mathcal{S}_{\text{IHT}}$. Of course, if we know $\beta$ and $H$, then all we need to do is to check at each step $j$ whether the smallest singular value of the matrix $(\beta, \ldots, H^j\beta)$, $j = 0, \ldots, p-1$, is zero, and stop as soon as it is. At that point, $k = j$. However, in practice, $H$ and $\beta$ are replaced by their sample estimates, say $\widehat{H}$ and $\hat{\beta}$. Whereas in the population sequence $\{\beta, H\beta, \ldots\}$ the smallest singular value becomes zero after a certain point, in the sample sequence $\{\hat{\beta}, \widehat{H}\hat{\beta}, \ldots\}$ the smallest singular value becomes small, rather than zero, after a certain point. So our task is to deduce an asymptotic distribution against which we can judge if the observed smallest singular values correspond to singular values of 0 in the population.

Procedures for determining the order of a dimension reduction space via sequential testing of hypotheses were developed previously for other dimension reduction methods. For example, Li (1992) developed a testing procedure for PHD, and Li (1991) and Schott (1994) developed testing procedures for SIR.

It should be mentioned that, although IHT can bring out multiple vectors in the CMS, at the present stage we do not have a rigorous set of sufficient conditions that guarantees IHT will actually cover the CMS. However, based on our numerous experiences with real data, IHT often does well in bringing out the full pattern in the conditional mean. Hence in part of the subsequent development (i.e., the constrained case) we take the pragmatic approach of making coverage a working assumption at the outset. Indeed, the issue of coverage is challenging, and to date there has not been a general result published in this regard. For this reason a similar working assumption is typically adopted for the asymptotic development of other methods, such as those for SIR and PHD [Li (1991, 1992)]. Or, alternatively, the null hypothesis is formulated directly on the rank of the population matrix corresponding to the estimator rather than on the dimension of the CS [Schott (1994)]. In fact, we are inclined to believe that IHT is more comprehensive in estimating the CMS than using OLS or PHD alone, because it can pick up both monotone and U-shaped trends, so long as it has a nonzero vector such as OLS to prime the process. Cook and Li (2002) argued that the span of OLS and PHD is a subset of the CMS, and it seems that a combination of them should provide a reasonably comprehensive estimator of the CMS. IHT can be viewed as one way of combining these elements, without evoking the constant covariance condition that is required by PHD.

The rest of the paper is organized as follows. In Section 2 we introduce a version of IHT that is modified slightly from that of Cook and Li (2002) for an invariance consideration. We also formulate the hypothesis testing problem and establish initial asymptotic expansions. In Section 3 we derive the asymptotic distribution under certain assumptions on the predictor $X$ and its relation with the response $Y$. In this case the limiting distribution is a chi-squared distribution (Theorem 4). In Section 4 we derive the asymptotic distribution when no practically restrictive conditions are placed on $X$ or $Y$ (Theorem 5). In Section 5 we discuss the implementation of the tests in both cases. We check our theoretical conclusions against simulated results in Section 6, and we apply both procedures to analyze an ozone data set.

## 2. Foundations.

2.1. *Invariant IHT.* Let $(X_1, Y_1), \ldots, (X_n, Y_n)$ be $n$ independent copies of $(X, Y)$, in which the predictor $X$ is a random vector in $\mathbb{R}^p$ and the response $Y$ is a scalar random variable. We assume throughout this article that $\mathrm{Var}(X)$ is positive definite. As demonstrated in Cook and Li (2002), the CMS is invariant under affine transformation of $X$: For any nonsingular $p$ by $p$ matrix $A$ and $p$-dimensional vector $b$,

$$\mathcal{S}_{E(Y|A^T X + b)} = A^{-1} \mathcal{S}_{E(Y|X)}.$$

Thus, without loss of generality, we prestandardize and use $Z = \mathrm{Var}(X)^{-1/2}(X - E(X))$ as the predictor vector so that $E(Z) = 0$ and $\mathrm{Var}(Z) = I_p$. In what follows, we first estimate the standardized subspace $\mathcal{S}_{E(Y|Z)}$ along with its dimension, and then transform back to estimate $\mathcal{S}_{E(Y|X)}$.

Let $\overline{Z}$ and $\widehat{\Sigma}$ be the sample mean and sample covariance matrix of $Z$:

$$\overline{Z} = E_n(Z) \quad \text{and} \quad \widehat{\Sigma} = E_n(Z - \overline{Z})(Z - \overline{Z})^T,$$

where $E_n f(Z)$ stands for $n^{-1} \sum_{i=1}^n f(Z_i)$. Let $\widehat{Z}$ be the standardized $Z$,

$$\widehat{Z} = \widehat{\Sigma}^{-1/2}(Z - \overline{Z}),$$

and let

$$\widetilde{H} = E_n\{\tilde{e} \widehat{Z} \widehat{Z}^T\},$$

where $\tilde{e}$ is the observed regression error $Y - \overline{Y} - \tilde{\beta}^T \widehat{Z}$ with $\tilde{\beta} = E_n(\widehat{Z}(Y - \overline{Y}))$. This matrix was suggested in Cook and Li (2002) as the transformation matrix in the r-based IHT method. However, in practice, it is desirable to make IHT invariant under affine transformation of both $Z$ and $Y$; that is, conclusions drawn from $(Z, Y)$ should be identical to those drawn from $(AZ + b, cY + d)$, where $A$ is any $p$ by $p$ nonsingular matrix, $b$ is any $p$-dimensional vector, $c$ is a nonzero scalar and $d$ is any scalar. For this purpose we will replace $Y - \overline{Y}$ in the above transformation by its standardized version

$$\widehat{Y} = \hat{\sigma}^{-1}(Y - \overline{Y}),$$

where $\hat{\sigma}^2 = E_n(Y - \overline{Y})^2$, and use the transformation matrix

$$\widehat{H} = E_n\{\hat{e}\widehat{Z}\widehat{Z}^T\},$$

where $\hat{e} = \widehat{Y} - \hat{\beta}^T \widehat{Z}$ with $\hat{\beta}$ being the regression estimate $E_n(\widehat{Y}\widehat{Z})$. For consistency in the rest of this article, we now redefine $H$ and $\beta$ to be the population versions of $\widehat{H}$ and $\hat{\beta}$,

$$H = E([(Y - E(Y))/\sigma - \beta^T Z]ZZ^T),$$

where $\beta = E((Y - E(Y))Z)/\sigma$.

2.2. *Formulation of hypotheses.* Let

$$B = (\beta, H\beta, \ldots, H^{p-1}\beta) \quad \text{and} \quad \widehat{B} = (\hat{\beta}, \widehat{H}\hat{\beta}, \ldots, \widehat{H}^{p-1}\hat{\beta}).$$

We estimate the rank of $B$, which is equal to the dimension of $\mathcal{S}_{\text{IHT}}$ since $\mathcal{S}_{\text{IHT}} = \text{Span}(B)$, by conducting a series of hypothesis tests. Let $\lambda_1 \geq \lambda_2 \geq \cdots \geq \lambda_p$ be the eigenvalues of $BB^T$, and consider the sequence of tests

$$H_{0,j} : \lambda_{j+1} = \cdots = \lambda_p = 0, \qquad j = 0, 1, \ldots, p - 1.$$

The rank $k$ of $B$ is the smallest value of $j$ for which this hypothesis holds. Let $\hat{\lambda}_1 \geq \hat{\lambda}_2 \geq \cdots \geq \hat{\lambda}_p$ be the eigenvalues of $n\widehat{B}\widehat{B}^T$. We test $H_{0,j}$ using the statistic

$$T_j = C^{-1} \sum_{i=j+1}^{p} \hat{\lambda}_i,$$

where $C$ is a positive constant that depends on $j$ and will be determined later. Relatively large values of $T_j$ provide evidence against $H_{0,j}$. Tests of $H_{0,j}$ are used to estimate the rank $k$ of $B$ as follows: Beginning with $j = 0$, test $H_{0,0}$. If the hypothesis is rejected, increment $j$ by one and test again, stopping with the first nonsignificant result. The corresponding value of $j$ is the estimate $\hat{k}$ of $k$. Procedures of this form are fairly common for estimating the rank of a matrix; see, for example, Rao [(1965), page 472].

2.3. *Initial asymptotic equivalences.* In this section we characterize the components of $T_k$ in terms of their asymptotically equivalent variables, and provide expansions that will be useful when studying the distribution of $T_k$ in later sections.

First, consider the singular value decomposition of $B$:

$$(1) \qquad B = (\Gamma_1 \ \Gamma_0) \begin{pmatrix} D & 0 \\ 0 & 0 \end{pmatrix} \begin{pmatrix} \Psi_1^T \\ \Psi_0^T \end{pmatrix},$$

where $\Gamma = (\Gamma_1, \Gamma_0)$ and $\Psi = (\Psi_1, \Psi_0)$ are $p$ by $p$ orthonormal matrices, $D$ is a $k$ by $k$ diagonal matrix with positive diagonal elements, $\Gamma_1$ and $\Psi_1$ have dimension $p$ by $k$, and $\Gamma_0$ and $\Psi_0$ have dimension $p$ by $p-k$. It follows from Eaton and Tyler (1994) that the joint asymptotic distribution of the $p-k$ smallest singular values of the matrix $\sqrt{n}\widehat{B}$ is the same as that of the singular values of the matrix $\sqrt{n}\Gamma_0^T(\widehat{B}-B)\Psi_0$. Therefore the asymptotic distribution of $CT_k$ is the same as that of

$$n\{\text{vec}[\Gamma_0^T(\widehat{B}-B)\Psi_0]\}^T \text{vec}[\Gamma_0^T(\widehat{B}-B)\Psi_0],$$

where vec is the usual operator that maps a matrix to a vector by stacking its columns: if $A$ is a matrix with columns $a_1, \ldots, a_p$, then $\text{vec}(A) = (a_1^T \cdots a_p^T)^T$. Thus determining the asymptotic distribution of $T_k$ boils down to computing the asymptotic distribution of $\sqrt{n}\,\text{vec}[\Gamma_0^T(\widehat{B}-B)\Psi_0]$.

The estimate $\widehat{B}$ is a function of $\hat{\beta}$ and $\widehat{H}$, both of which are essentially (though not exactly) sums of independent and identically distributed (i.i.d.) random variables. So the key is to expand $\sqrt{n}\,\text{vec}[\Gamma_0^T(\widehat{B}-B)\Psi_0]$ as a function of sums of i.i.d. random variables.

First, expand $\widehat{H}^i\hat{\beta}$ so that the remainder is of the order $O_p(n^{-1})$, $i = 1, \ldots, p-1$. Starting with

$$(2) \qquad \widehat{H}^i\hat{\beta} = \{H + (\widehat{H} - H)\}^i\{\beta + (\hat{\beta} - \beta)\},$$

the term $\{H + (\widehat{H} - H)\}^i$ can be expanded as the sum of $2^i$ terms, each being of the form $G_1 \cdots G_i$, where the $G$'s can be either $H$ or $\widehat{H} - H$. However, those $G_1 \cdots G_i$ terms involving two or more $\widehat{H} - H$ are of the order $O_p(n^{-1})$ or smaller and can be dropped. For the terms involving only one $\widehat{H} - H$, the $i-1$ $H$'s appear either on the left, or right, or both sides, of $\widehat{H} - H$. In other words they can be expressed as $H^j(\widehat{H} - H)H^{i-1-j}$, where $j = 0, \ldots, i-1$. Hence we have the following expansion:

$$\{H + (\widehat{H} - H)\}^i = H^i + \sum_{j=0}^{i-1} H^j(\widehat{H} - H)H^{i-1-j} + O_p(n^{-1}).$$

Substitute this expansion into (2) to obtain

$$\widehat{H}^i\hat{\beta} - H^i\beta = H^i(\hat{\beta} - \beta) + \sum_{j=0}^{i-1} H^j(\widehat{H} - H)H^{i-1-j}\beta + O_p(n^{-1}),$$

(3)

$$i=1,\ldots,p-1.$$

We next further expand $\hat{\beta} - \beta$ and $\widehat{H} - H$ as functions of sums of i.i.d. random variables. This is given in the next lemma; its proof is provided in the Appendix.

LEMMA 1. *Under regularity conditions we have the following expansions:*

(4)
$$\hat{\beta} - \beta = E_n(ZY) - \beta - \tfrac{1}{2}E_n(ZZ^T - I_p)\beta$$
$$- \tfrac{1}{2}E_n(Y^2 - 1)\beta + O_p(n^{-1}),$$

(5)
$$\widehat{H} - H = E_n\{e(ZZ^T - I_p) - H\} - \tfrac{1}{2}E_n(ZZ^T - I_p)H$$
$$- \tfrac{1}{2}HE_n(ZZ^T - I_p) - \tfrac{1}{2}E_n(Y^2 - 1)H + O_p(n^{-1}).$$

Here, the "regularity conditions" refer to those under which the central limit theorem applies to averages of i.i.d. random variables, which in our case are guaranteed if $Z$ has finite fourth moments.

In the next section we derive the asymptotic distribution of $T_k$ under a set of constraints on $X$ and $Y$. These constraints are similar to those imposed on SIR and PHD to produce chi-squared asymptotic distributions [Li (1991, 1992), Cook (1998b) and Bura and Cook (2001)]. We refer to this case as the *constrained case*. We then derive the asymptotic distribution without these conditions. While the results for the general case can be applied to the constrained case, the latter takes advantage of the structures imposed and performs better if the constraints are satisfied. It also has a simple form of a chi-squared distribution which is easy to use.

## 3. Asymptotic distribution for constrained case.

3.1. *Constraints.* In the constrained case we assume that:

(C1) The span of the IHT vectors exhausts the CMS, $\mathcal{S}_{\text{IHT}} = \mathcal{S}_{E(Y|Z)}$.
(C2) The predictor $Z$ is normally distributed.
(C3) $E(e^2|Z) = E(e^2|P_{\mathcal{S}_{E(Y|Z)}}Z)$, where $e = Y - \beta^T Z$ is the population regression error.

These assumptions are similar in spirit to those imposed on the constrained cases of PHD and SIR. Under the linearity condition (A) alone, $\mathcal{S}_{\text{IHT}} \subseteq S_{E(Y|Z)}$. In condition (C1) we carry this a step further and assume equality. This implies, for example, that if $\dim(\mathcal{S}_{E(Y|Z)}) > 0$, then we must have $\beta \neq 0$. Condition (C3) says that $\mathcal{S}_{E(e^2|Z)} \subseteq \mathcal{S}_{E(Y|Z)}$. That is, $\mathcal{S}_{E(Y|Z)}$ must be a mean subspace for the regression of $e^2$ on $Z$. This means that any heteroscedasticity present in the residuals must depend only on directions in the CMS. Conditions (C2) and (C3) are used to force a simple chi-squared asymptotic distribution for $T_k$.

Because $Z$ is normal, both the linearity condition (A) and the constant covariance condition (B) hold and thus $\mathrm{Span}(\beta, H) \subseteq \mathcal{S}_{E(Y|Z)}$. Also, for any integer $j > 0$, $H^j \beta \in \mathrm{Span}(H)$, implying that $\mathcal{S}_{\mathrm{IHT}} \subseteq \mathrm{Span}(\beta, H)$. Hence, it follows from condition (C1) that

$$\mathcal{S}_{\mathrm{IHT}} = \mathrm{Span}(\beta, H) = \mathcal{S}_{E(Y|Z)}. \tag{6}$$

Because r-based PHD [Li (1992)] is designed to estimate $\mathrm{Span}(H)$, it follows that in the constrained case IHT combines r-based PHD with OLS. Cook (1998b) found that y-based PHD is not very effective at finding linear trends and that the best results in practice are often found by informally combining OLS with r-based PHD. IHT is the first formal method for combining OLS with r-based PHD, making use of $\beta$ to find linear trends in the mean function and $H$ to find curvature.

We next consider the asymptotic distribution of $T_k$ in the constrained case, picking up the general argument at the end of Section 2.3.

3.2. *Expansion of* $\sqrt{n}\,\mathrm{vec}[\Gamma_0^T(\widehat{B} - B)\Psi_0]$. From the singular value decomposition (1) we know that $\Gamma_0^T B = 0$. That is, the columns of $\Gamma_0$ are orthogonal to the columns of $B$, and by (6) they are also orthogonal to the columns of $(\beta, H)$ because $\mathcal{S}_{\mathrm{IHT}} = \mathrm{Span}(B)$. Thus, if we multiply both sides of (3) by the matrix $\Gamma_0^T$ from the left, all the terms that begin with an $H$ drop, and we have

$$\Gamma_0^T(\widehat{H}^i \hat{\beta} - H^i \beta) = \Gamma_0^T(\widehat{H} - H)H^{i-1}\beta + O_p(n^{-1}).$$

It follows that

$$\begin{aligned}
&\sqrt{n}\Gamma_0^T(\widehat{B} - B)\Psi_0 \\
(7) \quad &= \sqrt{n}\Gamma_0^T(\hat{\beta} - \beta, (\widehat{H} - H)\beta, \ldots, (\widehat{H} - H)H^{p-2}\beta)\Psi_0 + O_p(n^{-1/2}) \\
&\equiv \sqrt{n}\Gamma_0^T(\hat{\beta} - \beta, (\widehat{H} - H)B_0)\Psi_0 + O_p(n^{-1/2}),
\end{aligned}$$

where $B_0 = (\beta, \ldots, H^{p-2}\beta)$.

Observe that, in expansions (4) and (5), the terms

$$-\beta, \quad -E_n(Y^2 - 1)\beta, \quad H, \quad HE_n(ZZ^T - I_p)/2, \quad E_n(Y^2 - 1)H/2$$

vanish if we multiply them by $\Gamma_0^T$ from the left. Therefore, $\Gamma_0^T(\widehat{B} - B)\Psi_0$ reduces to

$$\Gamma_0^T(E_n(ZY) - E_n(W)\beta/2, (E_n(eW) - E_n(W)H/2)B_0)\Psi_0 + O_p(n^{-1}),$$

where $W$ stands for the matrix $ZZ^T - I_p$. Using the relation $B = (\beta, HB_0)$ we can rewrite the above matrix as

$$\Gamma_0^T(E_n(ZY), E_n(eW)B_0)\Psi_0 - \tfrac{1}{2}\Gamma_0^T E_n(W)B\Psi_0 + O_p(n^{-1}).$$

Note that the second term drops because $B\Psi_0 = 0$. Furthermore, the identity matrix $I_p$ in $W = ZZ^T - I_p$ also drops because it is to be multiplied from the left by $\Gamma_0^T$ and from the right by $B_0$, and the columns of $\Gamma_0$ are orthogonal to the columns of $B_0$, which consists of the first $p-1$ columns of the matrix $B$. To conclude, we have the following expansion for $\sqrt{n}\Gamma_0^T(\widehat{B} - B)\Psi_0$.

THEOREM 1.  *Under conditions* (C1) *and* (C2) *for the constrained case, we have*

$$(8) \quad \sqrt{n}\Gamma_0^T(\widehat{B} - B)\Psi_0 = \sqrt{n}\Gamma_0^T(E_n(ZY), E_n(eZZ^T)B_0)\Psi_0 + O_p(n^{-1/2}).$$

The right-hand side of (8) can be further simplified using the properties of $\Psi_0$. We will do this in two separate cases, $\beta \in \mathrm{Span}(H)$ or $\beta \notin \mathrm{Span}(H)$, and then synthesize them into a simple and general formula.

3.3. *Case* I: $\beta \notin \mathrm{Span}(H)$.   The next lemma describes the structure of $\Psi_0$ in this case, which is the key to the simplification. Its proof is given in the Appendix.

LEMMA 2.   *If $\beta \notin \mathrm{Span}(H)$, then:*
(i) *the first row of $\Psi_0$ is a zero vector, and*
(ii) *the second row of $\Psi_0$ is not a zero vector.*

To simplify the expansion of $\sqrt{n}\Gamma_0^T(\widehat{B} - B)\Psi_0$ using this result, rewrite the expansion (8) as

$$(9) \quad \begin{aligned} \Gamma_0^T(\widehat{B} - B)\Psi_0 &= (E_n(\Gamma_0^T ZY), E_n(e\Gamma_0^T ZZ^T B_0))\Psi_0 + O_p(n^{-1}) \\ &= E_n(e\Gamma_0^T ZZ^T B_0)\Phi_0 + O_p(n^{-1}), \end{aligned}$$

where $\Phi_0$ is the $p-1$ by $p-k$ matrix comprising the second through the $p$th rows of the $p$ by $p-k$ matrix $\Psi_0$. Furthermore, because the first row of $\Psi_0$ is 0 and because $B\Psi_0 = 0$, we have

$$(H\beta, \ldots, H^{p-1}\beta)\Phi_0 = HB_0\Phi_0 = 0.$$

In other words, the columns of the matrix $B_0\Phi_0$ are orthogonal to the columns of $H$. Consequently, letting $Q_H = I_p - P_{\mathrm{Span}(H)}$,

$$(10) \quad B_0\Phi_0 = Q_H B_0\Phi_0 = (Q_H\beta, 0, \ldots, 0)\Phi_0 = Q_H\beta\alpha_0^T,$$

where $\alpha_0^T$ is the first row of the matrix $\Phi_0$, which by Lemma 2 is a nonzero vector.

Substitute (10) into the right-hand side of (9) to obtain

$$\Gamma_0^T(\widehat{B} - B)\Psi_0 = E_n(e\Gamma_0^T Z Z^T Q_H \beta \alpha_0^T) \equiv E_n(eUV^T) + O_p(n^{-1}),$$

where $U = \Gamma_0^T Z$ and $V = \alpha_0 \beta^T Q_H Z$. Hence

$$\sqrt{n}\operatorname{vec}\Gamma_0^T(\widehat{B} - B)\Psi_0 = \sqrt{n}\operatorname{vec}(E_n(eUV^T)) + O_p(n^{-1/2})$$
$$= \sqrt{n}E_n(eV \otimes U) + O_p(n^{-1/2}).$$

Letting the columns of the $p$ by $k$ matrix $\gamma$ be an orthonormal basis for $\mathcal{S}_{\text{IHT}}$, we see that $V = \alpha_0(\beta^T Z - \beta^T \gamma \gamma^T Z)$ is measurable with respect to $\gamma^T Z$ because $\beta \in \mathcal{S}_{\text{IHT}}$. This implies that $E(eV \otimes U) = 0$, as can be seen from the following derivation:

$$E(eV \otimes U) = E(E(e|Z)V \otimes U)$$
$$= E(E(e|\gamma^T Z)V \otimes U)$$
$$= E(eE(V \otimes U|\gamma^T Z))$$
$$= E(eV \otimes E(U|\gamma^T Z))$$
$$= E(eV) \otimes E(U) = 0,$$

where, for the second equality we used the definition of the CMS, for the fourth we used the measurability of $V$ with respect to $\gamma^T Z$, and for the fifth we used the independence between $U$ and $\gamma^T Z$, which follows from the normality of $Z$ and the orthogonality between the columns of $\Gamma_0$ and the columns of $\gamma$.

Hence, by the central limit theorem, $\sqrt{n}\,E_n(eV \otimes U)$ is asymptotically normal with mean 0 and covariance matrix $E\{e^2(V \otimes U)(V \otimes U)^T\}$. We now simplify this covariance matrix:

(11)
$$\begin{aligned}E\{e^2(V \otimes U)(V \otimes U)^T\} &= E\{e^2(VV^T \otimes UU^T)\}\\ &= E\{E(e^2|\gamma^T Z)(VV^T \otimes UU^T)\}\\ &= E\{e^2 E(VV^T \otimes UU^T|\gamma^T Z)\}\\ &= E\{(e^2 VV^T) \otimes E(UU^T|\gamma^T Z)\}\\ &= E(e^2 VV^T) \otimes E(UU^T)\\ &= E(e^2 VV^T) \otimes I_{p-k}.\end{aligned}$$

For the second equality we have used the assumption $E(e^2|Z) = E(e^2|\gamma^T Z)$, and for the last equality we have used the fact that $E(UU^T) = \Gamma_0^T E(ZZ^T)\Gamma_0 = I_{p-k}$. The rest of the equalities follow from the similar argument we used in the demonstration of $E(eU \otimes V) = 0$.

Now substitute the definition $V = \alpha_0 \beta^T Q_H Z$ into the expression $E(e^2 VV^T) \otimes I_{p-k}$:

$$(\beta^T Q_H E(e^2 ZZ^T) Q_H \beta)(\alpha_0 \alpha_0^T \otimes I_{p-k}) = E(e\beta^T Q_H Z)^2 (\alpha_0 \alpha_0^T \otimes I_{p-k})$$
$$\equiv C(\alpha_0 \alpha_0^T / \|\alpha_0\|^2) \otimes I_{p-k},$$

where $C = E(e\beta^T Q_H Z)^2 \|\alpha_0\|^2$. It is easy to verify that any matrix of the form $\alpha \alpha^T \otimes I_m$, where $\alpha$ is a unit vector, is an idempotent matrix of rank $m$. Therefore,

$$\sqrt{n}\Gamma_0^T (\widehat{B} - B)\Psi_0 / \sqrt{C} \xrightarrow{\mathcal{L}} N(0, R),$$

where $R$ is an idempotent matrix of rank $p - k$. So we have proved the following theorem.

THEOREM 2. *Let $\hat{\lambda}_1 \geq \cdots \geq \hat{\lambda}_p$ be the eigenvalues of the matrix $n\widehat{B}\widehat{B}^T$. Suppose that conditions* (C1)–(C3) *hold and that $\beta \notin \mathrm{Span}(H)$. Then, under the null hypothesis $H_{0,k}: \lambda_{k+1} = \cdots = \lambda_p = 0$, we have*

$$C^{-1} \sum_{i=k+1}^{p} \hat{\lambda}_i \xrightarrow{\mathcal{L}} \chi^2_{p-k},$$

*where $C = E(e\beta^T Q_H Z)^2 \|\alpha_0\|^2$, $\alpha_0$ being the second row of the matrix $\Psi_0$.*

3.4. *Case* II: $\beta \in \mathrm{Span}(H)$. As in the previous case, the special structure of the matrix $\Psi_0$ under $\beta \in \mathrm{Span}(H)$ plays a critical role in the simplification of the expansion of $\sqrt{n}\Gamma_0^T (\widehat{B} - B)\Psi_0$. This structure is described in the next lemma, and proved in the Appendix.

LEMMA 3. *If $\beta \in \mathrm{Span}(H)$, then the first row of $\Psi_0$ is not 0.*

Now write

$$\Psi_0 = \begin{pmatrix} \tau_0^T \\ \Phi_0 \end{pmatrix},$$

where $\tau_0$ is a vector in $\mathbb{R}^{p-k}$ and $\Phi_0$, as before, is a $p - 1$ by $p - k$ matrix. Since $B\Psi_0 = 0$, we have

$$HB_0 \Phi_0 + \beta \tau_0^T = 0.$$

Because $\beta \in \mathrm{Span}(H)$, $\beta = H\eta$ for some $\eta$ in $\mathbb{R}^p$. Hence

$$HB_0 \Phi_0 + H\eta \tau_0^T = H(B_0 \Phi_0 + \eta \tau_0^T) = 0.$$

That is, the columns of the matrix $B_0 \Phi_0 + \eta \tau_0^T$ are orthogonal to the rows (and hence columns) of $H$. Consequently,

$$B_0 \Phi_0 + \eta \tau_0^T = Q_H (B_0 \Phi_0 + \eta \tau_0^T),$$

where, as before, $Q_H$ is the projection matrix onto the orthogonal complement of $\mathrm{Span}(H)$. Since $\beta \in \mathrm{Span}(H)$, $\mathrm{Span}(B_0) \subseteq \mathrm{Span}(H)$ and therefore $Q_H B_0 \Phi_0 = 0$. Hence, letting "†" denote the Moore–Penrose generalized inverse,

$$B_0 \Phi_0 = -\eta \tau_0^T + Q_H \eta \tau_0^T = -(I - Q_H) \eta \tau_0^T$$
$$= -H(HH)^\dagger H \eta \tau_0^T = -H(HH)^\dagger \beta \tau_0^T.$$

From the definition of the Moore–Penrose generalized inverse of a symmetric matrix it is easy to see that $H(HH)^\dagger = H^\dagger$. Therefore,

(12) $$B_0 \Phi_0 = -H^\dagger \beta \tau_0^T.$$

Rewrite expansion (8) as

$$\sqrt{n}\Gamma_0^T(\widehat{B} - B)\Psi_0 = \sqrt{n}(E_n(YU), E_n(eUV^T))\Psi_0 + O_p(n^{-1/2}),$$

where $U$ is the $(p-k)$-dimensional vector $\Gamma_0^T Z$ and $V$ is the $(p-1)$-dimensional vector $B_0^T Z$. Note that the $V$ here is different from that defined in Section 3.3, but $U$ denotes the same quantity. Since the columns of $B_0$ belong to $\mathrm{Span}(B)$, $V$ is measurable with respect to $\gamma^T Z$, and since the columns of $\Gamma_0$ are orthogonal to $\mathrm{Span}(B)$, $U$ and $\gamma^T Z$ are independent. Thus, following the same argument used in Section 3.3 for the demonstration of $E(eV \otimes U) = 0$, we can show that $E(YU) = 0$ and $E(eUV^T) = 0$. Hence the vector

$$\sqrt{n}\,\mathrm{vec}\{(E_n(YU), E_n(eUV^T))\Psi_0\} = \sqrt{n}(\Psi_0^T \otimes I_{p-k})\begin{pmatrix} E_n(YU) \\ E_n(eV \otimes U) \end{pmatrix}$$

is asymptotically multivariate normal of dimension $(p-k)^2$ with mean 0 and variance matrix

(13) $(\Psi_0^T \otimes I_{p-k}) \begin{pmatrix} E(Y^2 UU^T) & E(eYU(V \otimes U)^T) \\ E(eY(V \otimes U)U^T) & E(e^2(V \otimes U)(V \otimes U)^T) \end{pmatrix} (\Psi_0 \otimes I_{p-k}).$

We next simplify this covariance matrix.

By an argument similar to (11), we can show that

$$E(Y^2 UU^T) = E(Y^2) I_{p-k},$$
$$E(e^2(V \otimes U)(V \otimes U)^T) = E(e^2 VV^T) \otimes I_{p-k}.$$

To derive $E(eYU(V \otimes U)^T)$, note that $U = 1 \otimes U$, where 1 is the scalar one. Hence

$$E(eYU(V \otimes U)^T) = E(eY(1 \otimes U)(V \otimes U)^T) = E(eY(V^T \otimes UU^T)).$$

Now apply the argument leading to (11) to obtain

$$E(eYU(V \otimes U)^T) = E(eYV^T) \otimes I_{p-k}.$$

Hence the asymptotic variance (13) now becomes

$$(\Psi_0^T A \Psi_0) \otimes I_{p-k} \qquad \text{where } A = \begin{pmatrix} E(Y^2) & E(eYV^T) \\ E(eYV) & E(e^2 VV^T) \end{pmatrix}.$$

Expressing $\Psi_0^T$ as $(\tau_0, \Phi_0^T)$, we can rewrite the matrix $\Psi_0^T A \Psi_0$ as

(14)
$$\begin{aligned} A &= E(Y^2)\tau_0 \tau_0^T + \tau_0 E(eYV^T)\Phi_0 \\ &\quad + \Phi_0^T E(eYV)\tau_0^T + \Phi_0^T E(e^2 VV^T)\Phi_0. \end{aligned}$$

Now recall that $B_0 \Phi_0 = -H^\dagger \beta \tau_0^T$. So

$$\Phi_0^T V = \Phi_0^T B_0^T Z = -\tau_0 \beta^T H^\dagger Z.$$

Substitute this relation into (14) to obtain

$$\begin{aligned} A &= \{E(Y^2) - 2E(eYZ^T)H^\dagger \beta + \beta^T H^\dagger E(e^2 ZZ^T) H^\dagger \beta\} \tau_0 \tau_0^T \\ &= E(Y - eZ^T H^\dagger \beta)^2 \, \tau_0 \tau_0^T \equiv C_1 \tau_0 \tau_0^T / \|\tau_0\|^2, \end{aligned}$$

where $C_1$ is the constant $E(Y - eZ^T H^\dagger \beta)^2 \|\tau_0\|^2$. Note that $\tau_0 \tau_0^T / \|\tau_0\|^2$ is an idempotent matrix of rank 1, and $(\tau_0 \tau_0^T / \|\tau_0\|^2) \otimes I_{p-k}$ is an idempotent matrix of rank $p-k$. Therefore

$$\sqrt{n} \Gamma_0^T (\widehat{B} - B) \Psi_0 / \sqrt{C_1} \xrightarrow{\mathcal{L}} N(0, R),$$

where $R$ is an idempotent matrix of rank $p-k$. We have proved the following theorem.

THEOREM 3. *Suppose that conditions* (C1)–(C3) *hold and that $\beta$ belongs to* Span($H$). *Then, under the null hypothesis* $H_{0,k} : \lambda_{k+1} = \cdots = \lambda_p = 0$, *we have*

$$C_1^{-1} \sum_{i=k+1}^p \hat{\lambda}_i \xrightarrow{\mathcal{L}} \chi^2_{p-k},$$

*where $C_1 = E(Y - eZ^T H^\dagger \beta)^2 \|\tau_0\|^2$, with $\tau_0^T$ being the first row of the matrix $\Psi_0$.*

3.5. *Synthesis of the two cases.* To apply directly the asymptotic results developed in Sections 3.3 and 3.4, one must determine at the outset whether $\beta$ belongs to Span($H$), which would likely be problematic in practice. In this section we derive a general result that synthesizes the two cases. This enables us to apply the test in the constrained case without having to know whether $\beta$ belongs to Span($H$) ahead of time.

Recall from (10) and (12) that

$$B_0 \Phi_0 = \begin{cases} Q_H \beta \alpha_0^T, & \text{if } \beta \notin \text{Span}(H), \\ -H^\dagger \beta \tau_0^T, & \text{if } \beta \in \text{Span}(H). \end{cases}$$

Using this relation we can rewrite the constants $C$ and $C_1$ as

$$
\begin{aligned}
(15) \quad C &= \operatorname{tr}\{\alpha_0 E(e^2 \beta^T Q_H Z Z^T Q_H \beta)\alpha_0^T\} \\
&= \operatorname{tr}\{\Phi_0^T B_0^T E(e^2 Z Z^T) B_0 \Phi_0\}, \\
(16) \quad C_1 &= \operatorname{tr}\{\tau_0 E(Y - e\beta^T H^\dagger Z)^2 \tau_0^T\} \\
&= \operatorname{tr}\{E(Y^2)\tau_0 \tau_0^T + 2\tau_0 E(YeZ^T) B_0 \Phi_0 + \Phi_0^T B_0^T E(e^2 Z Z^T) B_0 \Phi_0\}.
\end{aligned}
$$

Now consider the matrix

$$
(17) \quad A = \Psi_0^T \begin{pmatrix} 1 & 0 \\ 0 & B_0^T \end{pmatrix} \begin{pmatrix} E(Y^2) & E(YeZ^T) \\ E(YeZ) & E(e^2 Z Z^T) \end{pmatrix} \begin{pmatrix} 1 & 0 \\ 0 & B_0 \end{pmatrix} \Psi_0.
$$

If $\beta \notin \operatorname{Span}(H)$, then the first row of $\Psi_0$ is 0 and the matrix reduces to that inside trace$(\cdot)$ on the right-hand side of (15). If $\beta \in \operatorname{Span}(H)$, then the first row of $\Psi_0$ is $\tau_0$ and the matrix reduces to that inside trace$(\cdot)$ of (16). Thus if we let $C_2$ be $\operatorname{tr}(A)$, then it automatically normalizes the asymptotic distribution of $T_k$ to a $\chi^2_{p-k}$ distribution in both cases. To further simplify the notation, let $W$ denote the $(p+1)$-dimensional vector $(Y, eZ^T)^T$, and let $\operatorname{diag}(1, B_0)$ denote the $p+1 \times p$ block-diagonal matrix with diagonal blocks 1 (the scalar one) and $B_0$. Then we can express $C_2$ as

$$
(18) \quad C_2 = \operatorname{tr}\{\Psi_0^T \operatorname{diag}(1, B_0^T) E(WW^T) \operatorname{diag}(1, B_0) \Psi_0\}.
$$

The next theorem summarizes this general result.

THEOREM 4. *Suppose that conditions* (C1)–(C3) *hold. Then, under the null hypothesis* $H_{0,k} : \lambda_{k+1} = \cdots = \lambda_p = 0$, *we have*

$$
C_2^{-1} \sum_{i=k+1}^{p} \hat{\lambda}_i \xrightarrow{\mathcal{L}} \chi^2_{p-k},
$$

*where* $C_2$ *is defined at* (18).

**4. Asymptotic distribution for the general case.** We now derive the asymptotic distribution of $C_2^{-1} \sum_{i=k+1}^{p} \hat{\lambda}_i$ in the general case. The asymptotic result of this section holds without conditions (C1)–(C3)—in its general form the asymptotic distribution is related only to the rank of the matrix $B$. Thus for clarity we will not refer to these conditions in the statement of the result. Under the linear conditional mean condition (A), $\mathcal{S}_{\text{IHT}}$ is a subspace of the CMS, and the test helps us to identify a set of significant vectors that belongs to the CMS. Under the coverage condition (C1), $\mathcal{S}_{\text{IHT}}$ is equal to the CMS, and the test helps us to identify the CMS itself. The point of this generalization is that (C3) is altogether removed, (C2) is replaced by the much weaker condition (A), and without (C1) we can still find vectors in the CMS but without the guarantee that they will span the CMS.

By expansion (3), the leading term of $\widehat{B} - B$, ignoring the error of magnitude $O_p(n^{-1})$, is

$$\left(\hat{\beta} - \beta, H(\hat{\beta} - \beta) + (\widehat{H} - H)\beta, \ldots, H^{p-1}(\hat{\beta} - \beta)\right.$$
$$\left. + \sum_{j=0}^{p-2} H^j(\widehat{H} - H)H^{p-2-j}\beta\right).$$

This can be written as the sum of $p$ matrices of simpler structures, as follows:

$$(\hat{\beta} - \beta, \ldots, H^{p-1}(\hat{\beta} - \beta)) + (0, (\widehat{H} - H)\beta, \ldots, (\widehat{H} - H)H^{p-2}\beta)$$
$$+ (0, 0, H(\widehat{H} - H)\beta, \ldots, H(\widehat{H} - H)H^{p-3}\beta) + \cdots$$
$$+ (0, \ldots, 0, H^{p-2}(\widehat{H} - H)\beta).$$

Thus, the vector $\text{vec}(\widehat{B} - B)$ can be written as

$$\begin{pmatrix} \hat{\beta} - \beta \\ H(\hat{\beta} - \beta) \\ \vdots \\ H^{p-1}(\hat{\beta} - \beta) \end{pmatrix} + \begin{pmatrix} 0 \\ (\widehat{H} - H)\beta \\ \vdots \\ (\widehat{H} - H)H^{p-2}\beta \end{pmatrix} + \cdots + \begin{pmatrix} 0 \\ 0 \\ \vdots \\ H^{p-2}(\widehat{H} - H)\beta \end{pmatrix}.$$

In other words,

$$(19) \quad \text{vec}(\widehat{B} - B) = \begin{pmatrix} I_p & 0 & \cdots & 0 & 0 \\ H & I_p & \cdots & 0 & 0 \\ \vdots & & \ddots & & \vdots \\ H^{p-1} & & \cdots & H & I_p \end{pmatrix} \times \begin{pmatrix} \hat{\beta} - \beta \\ (\widehat{H} - H)\beta \\ \vdots \\ (\widehat{H} - H)H^{p-2}\beta \end{pmatrix}$$
$$\equiv MV.$$

Here, $M$ is a $p^2$ by $p^2$ constant matrix, and $V$ is a random vector consisting of subvectors $\hat{\beta} - \beta, \ldots, (\widehat{H} - H)H^{p-2}\beta$ which, according to Lemma 1, can be approximated by sums of independent and identically distributed vectors. It follows that $\sqrt{n}MV$ converges in distribution to a $p^2$-dimensional multivariate normal random vector.

To write an explicit form of the asymptotic distribution, let

$$\xi_1 = ZY - \beta - (ZZ^T - I_p)\beta/2 - (Y^2 - 1)\beta/2,$$
$$(20) \quad \xi_i = \{e(ZZ^T - I_p) - H - (ZZ^T - I_p)H/2$$
$$- H(ZZ^T - I_p)/2 - (Y^2 - 1)H/2\}H^{i-2}\beta, \quad i = 2, \ldots, p.$$

Then, by Lemma 1,
$$\hat{\beta} - \beta = E_n(\xi_1) + O_p(n^{-1}),$$
$$(\widehat{H} - H)H^{i-2}\beta = E_n(\xi_i) + O_p(n^{-1}), \qquad i = 2, \ldots, p.$$

Thus, if we let $\xi$ be the vector $(\xi_1^T, \ldots, \xi_p^T)^T$, then $V = E_n(\xi) + O_p(n^{-1})$, and consequently $\sqrt{n}V$ converges in distribution to a $p^2$-dimensional multivariate normal with mean 0 and covariance matrix $E(\xi\xi^T)$. Therefore,
$$\sqrt{n/C_2}\operatorname{vec}\{\Gamma_0^T(\widehat{B} - B)\Psi_0\} \xrightarrow{\mathcal{L}} N(0, (\Psi_0 \otimes \Gamma_0)^T M E(\xi\xi^T) M^T (\Psi_0 \otimes \Gamma_0)/C_2).$$

Thus we have proved the following theorem.

THEOREM 5. *In the general case, we have*
$$C_2^{-1} \sum_{i=k+1}^{p} \hat{\lambda}_i \xrightarrow{\mathcal{L}} \sum_{i=1}^{(p-k)^2} \omega_i K_i,$$
*where $K_1, \ldots, K_{(p-k)^2}$ are independent chi-squared random variables with one degree of freedom and $\omega, \ldots, \omega_{(p-k)^2}$ are the eigenvalues of the matrix*

(21) $$(\Psi_0 \otimes \Gamma_0)^T M E(\xi\xi^T) M^T (\Psi_0 \otimes \Gamma_0)/C_2,$$

*with $C_2$, $M$ and $\xi$ defined by* (18), (19) *and* (20), *respectively.*

**5. Implementation.** In this section we describe how to estimate the various unknown quantities involved in the asymptotic distribution of $T_k = C_2^{-1}\sum_{i=k+1}^{p}\hat{\lambda}_i$, for both the constrained and the general cases. In the constrained case we only need to estimate $C_2$, whereas in the general case we need also to estimate the coefficients $\omega_1, \ldots, \omega_{(p-k)^2}$.

To estimate $C_2$, recall that $\Gamma_0$ and $\Psi_0$ are derived from the singular value decomposition of $B$. That is, the columns of $\Gamma_0$ are the eigenvectors of the matrix $BB^T$ corresponding to its zero eigenvalues, and the columns of $\Psi_0$ are the eigenvectors of $B^T B$ corresponding to its zero eigenvalues. So, we let $\widehat{\Gamma}_0$ be the $p \times p - j$ matrix whose columns are the $p - j$ eigenvectors of $\widehat{B}\widehat{B}^T$ corresponding to its smallest eigenvalues, in a descending order. In a similar manner construct $\widehat{\Psi}_0$, also of dimension $p \times p - j$, from the matrix $\widehat{B}^T\widehat{B}$. Furthermore, we will estimate $B_0$ by its sample version
$$\widehat{B}_0 = (\hat{\beta}, \ldots, \widehat{H}^{p-2}\hat{\beta}).$$
By the weak law of large numbers $\widehat{B}$ and $\widehat{B}_0$ consistently estimate $B$ and $B_0$, and, under the null hypothesis $H_{0,j}$, the matrices $\widehat{\Gamma}_0$ and $\widehat{\Psi}_0$ consistently estimate $\Gamma_0$ and $\Psi_0$. We propose to estimate $C_2$ by substituting the estimates $\widehat{\Gamma}_0$, $\widehat{\Psi}_0$ and $\widehat{B}_0$ for their population values $\Gamma_0$, $\Psi_0$ and $B_0$ in (18):
$$\widehat{C}_2 = \operatorname{tr}\{\widehat{\Psi}_0^T \operatorname{diag}(1, \widehat{B}_0^T) E_n(\widehat{W}\widehat{W}^T) \operatorname{diag}(1, \widehat{B}_0)\widehat{\Psi}_0\},$$

where $\widehat{W}$ is the $(p+1)$-dimensional vector $(\widehat{Y}, \hat{e}\widehat{Z}^T)^T$. By Slutsky's theorem, substituting $\widehat{C}_2$ for $C_2$ in $T_j = C_2^{-1} \sum_{i=j+1}^{p} \hat{\lambda}_i$ will not change its asymptotic distribution.

To estimate $\omega_1, \ldots, \omega_{(p-k)^2}$, let

$$\hat{\xi}_1 = \widehat{Z}\widehat{Y} - \hat{\beta} - (\widehat{Z}\widehat{Z}^T - I_p)\hat{\beta}/2 - (\widehat{Y}^2 - 1)\hat{\beta}/2,$$
$$\hat{\xi}_i = \{\hat{e}(\widehat{Z}\widehat{Z}^T - I_p) - \widehat{H} - (\widehat{Z}\widehat{Z}^T - I_p)/2$$
$$\quad - (\widehat{Z}\widehat{Z}^T - I_p)\widehat{H}/2 - (\widehat{Y}^2 - 1)\widehat{H}/2\}\widehat{H}^{i-2}\hat{\beta}, \qquad i = 2, \ldots, p,$$

and let $\hat{\xi}$ be the $p^2$-dimensional vector $(\hat{\xi}_1^T, \ldots, \hat{\xi}_p^T)^T$. We estimate $M$ in (21) by replacing $H$ in the definition of $M$, as given in (19), by $\widehat{H}$. The coefficients $\omega_1, \ldots, \omega_{(p-k)^2}$ are then estimated by the eigenvalues of

$$(\widehat{\Psi}_0 \otimes \widehat{\Gamma}_0)^T \widehat{M} E_n(\hat{\xi}\hat{\xi}^T) \widehat{M}^T (\widehat{\Psi}_0 \otimes \widehat{\Gamma}_0)/\widehat{C}_2.$$

As we can see from its construction, the general test does not reduce numerically to the constrained case when conditions (C1)–(C3) are satisfied, though in this case the two asymptotic distributions are first-order equivalent because $\hat{\omega}_1, \ldots, \hat{\omega}_{(p-k)^2}$ converges in probability to $\omega_1, \ldots, \omega_{(p-k)^2}$, which contain $p - k$ 1's and $(p - k)^2 - (p - k)$ 0's. Because the test for the constrained case uses this special 0–1 structure of the $\omega$'s, it is expected to outperform the general test when conditions (C1)–(C3) are satisfied.

## 6. Simulation results.

6.1. *Test levels.* In this section we present selected results from a simulation study to investigate the levels of the chi-squared and weighted chi-squared tests, one goal being to provide support for the validity of our asymptotic results. The tests have the same statistic for the hypothesis

$$\dim(\mathcal{S}_{\mathrm{IHT}}) = j, \qquad T_j = \widehat{C}_2^{-1} \sum_{i=j+1}^{p} \hat{\lambda}_i,$$

but use different reference distributions. The $\chi^2_{p-j}$ reference distribution is appropriate under conditions (C1)–(C3) of Theorem 4. Otherwise the reference distribution is the weighted chi-squared of Theorem 5. The scaling constant $C_2$ and the weights $\omega_i$ for the weighted chi-squared reference distribution were estimated as indicated in Section 5. For each simulation run the estimated test levels were based on 1000 replications and where relevant the chi-squared and weighted chi-squared tests were performed on the same data. There is a substantial literature on computing tail areas of distributions of linear combinations of chi-squared random variables. See Field (1993) for an introduction. For reference, the nominal standard errors of the estimated levels of nominal 1, 5, 10 and 15 percent tests are about 0.31, 0.69, 0.95 and 1.13.

Table 1 contains results for a null regression with four independent standard normal predictors and an independent standard normal response. The estimated levels of the tests seem quite far from the nominal levels for $n = 25$ observations, but the agreement seems good for both tests with more than about $n = 100$ observations.

Tables 2 and 3 contain results based on the model

(22) $$Y = Z_1 + 0.2(Z_1 + Z_2)^2 + \sigma N(0, 1)$$

with $p = 4$ independent standard normal predictors, and various sample sizes and values for $\sigma$. Because $\dim(\mathcal{S}_{\mathrm{IHT}}) = 2$, we studied the behavior of $T_2$ in each case. In the first part of Table 2 we held the sample size fixed at $n = 50$ and varied $\sigma$ from 0 to 1.6. As $\sigma$ increases the estimated levels of both tests tend to decrease, ending with quite conservative tests at $\sigma = 1.6$. Note, however, that $\sigma = 1.6$ is large compared to the variance of $Z_1$, with $\mathrm{var}(\varepsilon)/\mathrm{var}(Z_1) = 2.56$. At this error rate it is not surprising that the percentiles for both tests differ quite a bit from their nominal value, because the power of $T_1$ is not much larger than the nominal error rate. In the second part of Table 2 we held $\sigma$ fixed at 1.6 and increased the sample size. As the sample size increases we see the asymptotic approximations improving, ending with reasonable results at $n = 400$. Our general conclusion from Tables 1 and 2 and other simulation results not reported here is that the results are behaving as expected, which supports our analytic calculations and method of implementation suggested in Section 5. Perhaps as expected, the weighted chi-squared seems to take a larger sample size for the asymptotics to take

TABLE 1
*Estimated level of nominal* 1, 5, 10 *and* 15 *percent chi-squared* ($\chi^2$) *and weighted chi-squared* ($\bar{\chi}^2$) *tests based on* $T_0$ *for a* 0D *regression with* $p = 4$ *independent standard normal predictors and an independent standard normal response*

|         |               | Nominal level (%) |     |      |      |
|---------|---------------|-------------------|-----|------|------|
| $n$     | Test          | 1                 | 5   | 10   | 15   |
| 25      | $\chi^2$      | 0                 | 2.4 | 8.2  | 15.3 |
| 25      | $\bar{\chi}^2$| 2.4               | 8.7 | 15.8 | 21.6 |
| 50      | $\chi^2$      | 0.1               | 2.4 | 9.4  | 15.4 |
| 50      | $\bar{\chi}^2$| 0.8               | 6.6 | 12.1 | 18.2 |
| 100     | $\chi^2$      | 0.8               | 4.6 | 9.9  | 15.4 |
| 100     | $\bar{\chi}^2$| 1.5               | 6.4 | 11.8 | 17.4 |
| 200     | $\chi^2$      | 1.1               | 4.1 | 9.3  | 14.5 |
| 200     | $\bar{\chi}^2$| 1.3               | 4.3 | 10.5 | 14.8 |

Table 2

*Estimated level of nominal* 1, 5, 10 *and* 15 *percent chi-squared* ($\chi^2$) *and weighted chi-squared* ($\bar{\chi}^2$) *tests based on* $T_2$ *for model* (22) *with* $p = 4$ *independent standard normal predictors* $Z_j$

|          |                | Nominal level, $n = 50$ |     |      |      |
|----------|----------------|------|-----|------|------|
| $\sigma$ | Test           | 1    | 5   | 10   | 15   |
| 0        | $\chi^2$       | 1.1  | 4.8 | 9.5  | 14.6 |
| 0        | $\bar{\chi}^2$ | 0.1  | 3.3 | 8.3  | 13.7 |
| 0.2      | $\chi^2$       | 1.0  | 5.9 | 10.4 | 15.1 |
| 0.2      | $\bar{\chi}^2$ | 0.3  | 3.8 | 8.7  | 15.4 |
| 0.4      | $\chi^2$       | 1.1  | 4.7 | 10.0 | 14.6 |
| 0.4      | $\bar{\chi}^2$ | 0    | 1.9 | 5.7  | 10.6 |
| 0.8      | $\chi^2$       | 0.5  | 3.5 | 7.7  | 11.4 |
| 0.8      | $\bar{\chi}^2$ | 0    | 1.0 | 3.7  | 6.3  |
| 1.6      | $\chi^2$       | 0    | 0.7 | 3.1  | 5.7  |
| 1.6      | $\bar{\chi}^2$ | 0    | 0.2 | 0.7  | 1.6  |

|       |                | Nominal level, $\sigma = 1.6$ |     |      |      |
|-------|----------------|------|-----|------|------|
| $n$   | Test           | 1    | 5   | 10   | 15   |
| 100   | $\chi^2$       | 0.1  | 2.9 | 6.3  | 9.4  |
| 100   | $\bar{\chi}^2$ | 0.1  | 0.9 | 1.9  | 4.1  |
| 200   | $\chi^2$       | 1.3  | 5.7 | 10.4 | 14.2 |
| 200   | $\bar{\chi}^2$ | 0    | 1.9 | 5.0  | 8.6  |
| 400   | $\chi^2$       | 1.4  | 4.6 | 9.4  | 15.5 |
| 400   | $\bar{\chi}^2$ | 0.7  | 4.1 | 9.5  | 14.2 |

hold. Additionally, there is a tendency for the weighted chi-squared test to be conservative.

Table 3

*Estimated level of nominal* 1, 5, 10 *and* 15 *percent chi-squared tests based on* $T_2$ *for model* (22) *with* $p$ *independent standard normal predictors* $Z_j$

|     | Nominal level, $n = 100$ |     |      |      |
|-----|------|-----|------|------|
| $p$ | 1    | 5   | 10   | 15   |
| 4   | 1.1  | 4.7 | 10.0 | 14.2 |
| 6   | 0.7  | 4.9 | 9.9  | 15.0 |
| 8   | 0.9  | 4.0 | 8.6  | 15.9 |
| 12  | 0.3  | 4.6 | 10.2 | 15.1 |
| 16  | 0.5  | 3.9 | 7.6  | 12.0 |

In Table 3 we investigate the impact on the chi-squared test of increasing the number of unimportant predictors, holding $n = 100$ and $\sigma = 0.2$. Although there seems to be a little tendency for the estimated levels to decrease as $p$ increases, overall increasing the number of predictors does not seem to have much of an impact.

In Table 4 we consider model (22) with the error $\sigma N(0,1)$ term replaced by $0.5(\chi_2^2 - 2)$. Replacing the normal error with a chi-squared error did not seem to have a notable impact on the results. Because the error does not satisfy (C2) we have used the weighted chi-squared reference distribution.

Finally, we present a few confirmatory results based on $p = 5$ standard normal predictors and a response generated as

$$(23) \qquad Y = e^{0.3(2Z_1 + 3Z_2)} + 1.6\sin(Z_1 - Z_2) + \sigma N(0,1).$$

Letting SD denote a population standard deviation, the signal-to-noise ratio $\mathrm{SD}(\mathrm{E}(Y|Z))/\sigma$ for model (23) is about 0.4 times that for model (22), so the mean function of (23) should be harder to estimate. Table 5 contains estimated levels of $T_2$ for model (23) for three sample sizes and four

TABLE 4
*Estimated level of nominal* 1, 5, 10 *and* 15 *percent weighted chi-squared tests based on* $T_2$ *for a* 2D *regression with* $p = 4$ *independent standard normal predictors* $Z_j$ *and response* $Y = Z_1 + 0.2(Z_1 + Z_2)^2 + 0.5(\chi_2^2 - 2)$

|       | Nominal level (%) |     |     |      |
|-------|-------|-----|-----|------|
| $n$   | 1     | 5   | 10  | 15   |
| 50    | 0.1   | 2.5 | 6.5 | 11.0 |
| 100   | 0.5   | 3.8 | 8.9 | 13.1 |
| 200   | 1.0   | 5.8 | 9.9 | 14.6 |

TABLE 5
*Estimated level of nominal* 1, 5, 10 *and* 15 *percent chi-squared* $(\chi^2)$ *and weighted chi-squared* $(\bar{\chi}^2)$ *tests based on* $T_2$ *for simulation model* (23) *with* $\sigma = 0.2$

|     |            | Nominal level (%) |     |      |      |
|-----|------------|-----|-----|------|------|
| $n$ | Test       | 1   | 5   | 10   | 15   |
| 50  | $\chi^2$   | 0.5 | 4.1 | 9.5  | 13.7 |
| 50  | $\bar{\chi}^2$ | 0 | 2.1 | 4.6  | 8.7  |
| 100 | $\chi^2$   | 0.6 | 3.3 | 7.2  | 12.1 |
| 100 | $\bar{\chi}^2$ | 0.2 | 1.2 | 4.5 | 8.2  |
| 200 | $\chi^2$   | 1.2 | 4.5 | 10.5 | 14.5 |
| 200 | $\bar{\chi}^2$ | 0.4 | 2.4 | 6.3 | 11.5 |

TABLE 6
*Distribution of dimension estimates $\hat{k}$ out of 1000 trials based on the sequential chi-squared $(\chi^2)$ and weighted chi-squared $(\bar{\chi}^2)$ tests with constant nominal level $\alpha$*

|          | $\hat{k}$ with $n=50$ |     |     |        | $\hat{k}$ with $n=100$ |     |     |        |
| -------- | --- | --- | --- | ------ | --- | --- | --- | ------ |
| $\alpha$ | 0   | 1   | 2   | $\geq 3$ | 0   | 1   | 2   | $\geq 3$ |
| $\chi^2$ |     |     |     |        |     |     |     |        |
| 0.001    | 4   | 924 | 72  | 0      | 0   | 322 | 678 | 0      |
| 0.01     | 0   | 631 | 365 | 4      | 0   | 85  | 904 | 11     |
| 0.05     | 0   | 308 | 656 | 36     | 0   | 20  | 931 | 49     |
| 0.10     | 0   | 187 | 731 | 82     | 0   | 9   | 901 | 90     |
| 0.15     | 0   | 136 | 733 | 131    | 0   | 3   | 859 | 138    |
| $\bar{\chi}^2$ |  |  |  |  |  |  |  |  |
| 0.001    | 171 | 692 | 137 | 0      | 81  | 412 | 507 | 0      |
| 0.01     | 67  | 527 | 404 | 2      | 37  | 149 | 809 | 5      |
| 0.05     | 18  | 284 | 666 | 32     | 9   | 38  | 911 | 42     |
| 0.10     | 5   | 167 | 746 | 82     | 4   | 18  | 899 | 79     |
| 0.15     | 2   | 119 | 755 | 124    | 1   | 7   | 867 | 125    |

The model is (22) with $\sigma = 0.4$ and sample sizes 50 and 100.

nominal levels. These results are qualitatively similar to those discussed previously and confirm the conservative nature of the weighted chi-squared test in smaller samples.

6.2. *Estimation of* $\dim(\mathcal{S}_{\text{IHT}})$. In this section we present first results on the behavior of the sequential testing procedure discussed in Section 2.2 for estimating $k = \dim(\mathcal{S}_{\text{IHT}})$. We consider only estimates based on using the same nominal level for each of the sequential tests, although in a more comprehensive investigation it might be desirable to include variable levels.

Reasoning in the context of model (22) with $\dim(\mathcal{S}_{\text{IHT}}) = 2$, if the leading tests of $k = 0$ and $k = 1$ have power 1, then all of the estimation error arises from the level $\alpha$ of the test of $k = 2$, resulting in estimates $\hat{k} = 2$ with probability $1 - \alpha$ and $\hat{k} > 2$ with probability $\alpha$. Ideally, we would like to make $\alpha$ small, while maintaining high power in the leading tests. Leading tests with small values of $\alpha$ will have relatively low power and will tend to result in underestimation of $k$. We can increase the power of the leading tests by increasing $\alpha$, but this also increases the probability of overestimation.

For instance, Table 6 gives the empirical distribution of $\hat{k}$ out of 1000 trials based on the sequential chi-squared $(\chi^2)$ and weighted chi-squared $(\bar{\chi}^2)$ tests for model (22). With $n = 50$ we would prefer a level around 0.1 since the fraction of correct decisions was observed to change little for $\alpha > 0.1$ until it began to decrease. With $n = 100$, a level around 0.05 tends to balance over- and underestimation and produce the best results. With larger sample sizes or smaller values of $\sigma$, a level less than 0.05 may be preferred. Results for model (23) were qualitatively similar, but not quite as strong since its mean function is harder to estimate.

Overall, we found no compelling reason to prefer estimates with $\alpha > 0.15$. Tests with $\alpha = 0.05$ or $\alpha = 0.1$ tended to produce good results in our simulations, but tests with $\alpha < 0.05$ might yield better estimates with a significantly larger sample size or stronger signal.

6.3. *Direction estimation.* Given $k = 2$ for models (22) and (23), we studied the accuracy of the IHT estimates of the CMS by computing the absolute correlation between $Z_j$ and the fitted values from the OLS regression of $Z_j$ on the first two IHT predictors, $j = 1, 2$. Shown in Table 7 are three quantiles of the empirical distribution of these correlations over 1000 simulations for model (22). The results for model (23) are qualitatively similar, but as expected the correlations are smaller at the same sample size and standard deviation $\sigma$. For example, the quantiles for $Z_1$ under model (23) with $n = 100$ and $\sigma = 0.4$ were observed to be $q_{0.05} = 0.89$, $q_{0.5} = 0.97$ and $q_{0.95} = 0.996$.

**7. Ozone data.** In addition to simulations with normal predictors, we analyzed several different simulated and real data sets with nonnormal predictors using IHT methodology and found that in nearly all cases the chi-squared and weighted chi-squared reference distributions result in the same estimate of the dimension of the CMS. The ozone data [Breiman and Friedman (1985)] considered briefly in this section is an instance where the estimates of dimension differ.

TABLE 7
*Quantiles ($q_{0.05}$, $q_{0.5}$ and $q_{0.95}$) of the empirical distribution of the absolute correlation between $Z_j$ and the fitted values from the OLS regression of $Z_j$ on the first two IHT predictors, $j = 1, 2$, based on 1000 replications from model* (22) *with three values for $\sigma$ and two sample sizes*

| $\sigma$ | $n = 50$ | | | $n = 100$ | | |
|---|---|---|---|---|---|---|
| | $q_{0.05}$ | $q_{0.5}$ | $q_{0.95}$ | $q_{0.05}$ | $q_{0.5}$ | $q_{0.95}$ |
| | | | (A) $Z_1$ | | | |
| 0.2 | 0.98 | 0.995 | 0.9996 | 0.99 | 0.998 | 0.9998 |
| 0.4 | 0.97 | 0.993 | 0.9994 | 0.98 | 0.997 | 0.9997 |
| 0.8 | 0.94 | 0.99 | 0.999 | 0.97 | 0.994 | 0.9996 |
| | | | (B) $Z_2$ | | | |
| 0.2 | 0.77 | 0.95 | 0.996 | 0.89 | 0.98 | 0.998 |
| 0.4 | 0.71 | 0.94 | 0.996 | 0.86 | 0.97 | 0.997 |
| 0.8 | 0.47 | 0.88 | 0.992 | 0.72 | 0.94 | 0.995 |

The response $Y$ is atmospheric ozone concentration, and there are seven predictors: Daggett pressure gradient (DGPG, mmHg), humidity (HMDT, percent), visibility (VSTY, miles), wind speed (WDSP, mph), Vandenburg 500 millibar height (VDHT, m), the logarithms of Sandburg Air Force Base temperature (SBTP, degrees C), and inversion base temperature (IBTP, degrees F). The logarithm of the two temperature predictors was used to help to ensure that the linearity condition (A) hold to a useful approximation. Because we use only IHT methodology, there is no reason to consider the constant covariance condition (B).

The test results from using the chi-squared and weighted chi-squared reference distributions for $T_j$ are shown in Table 8. Use of the chi-squared reference distribution indicates that the dimension of the CMS is 3, while the weighted chi-squared reference distribution indicates two dimensions.

Shown in Figure 1 is a scatterplot of the response versus the first IHT predictor $\hat{v}_1^T \widehat{Z}$, where $\hat{v}_j$ is the eigenvector of $n\widehat{B}\widehat{B}^T$ corresponding to its $j$th largest eigenvalue $\hat{\lambda}_j$. A 3D plot (not shown) of the residuals $\hat{e}$ versus the first two IHT predictors $(\hat{v}_1^T \widehat{Z}, \hat{v}_2^T \widehat{Z})$ exhibits a saddle, confirming that $\dim(\mathcal{S}_{\text{IHT}})$ is at least 2. We were unable to find any notable graphical support for a third IHT predictor and consequently we conjecture that the results of the third chi-squared test in Table 8 are due to a failure of condition (C2) or (C3). In any event, because the conditions needed for the weighted chi-squared reference distribution are considerably less restrictive than those needed for the chi-squared, the weighted chi-squared $p$-values are likely more reliable.

**8. Discussion.** In this article we developed two asymptotic tests for the dimension $k$ of the IHT subspace $\mathcal{S}_{\text{IHT}}$. The tests use the same statistic $T_j = \widehat{C}_2^{-1} \sum_{i=j+1}^{p} \hat{\lambda}_i$ for the hypothesis $\text{rank}(B) = j$, but have different reference distributions depending on characteristics of the regression. The $\chi^2_{p-j}$ reference distribution is appropriate under conditions (C1)–(C3) of Theorem 4. Otherwise, in practically full generality, the reference distribution is the weighted chi-squared of Theorem 5. We typically use both reference distributions in practice, as illustrated in Table 8.

TABLE 8
*Test results for the ozone data*

| $j$ | $T_j$ | df | $\chi^2$ $p$-value | $\bar{\chi}^2$ $p$-value |
|---|---|---|---|---|
| 0 | 179.0 | 7 | 0.000 | 0.000 |
| 1 | 19.08 | 6 | 0.004 | 0.025 |
| 2 | 12.52 | 5 | 0.028 | 0.261 |
| 3 | 2.238 | 4 | 0.692 | 0.721 |

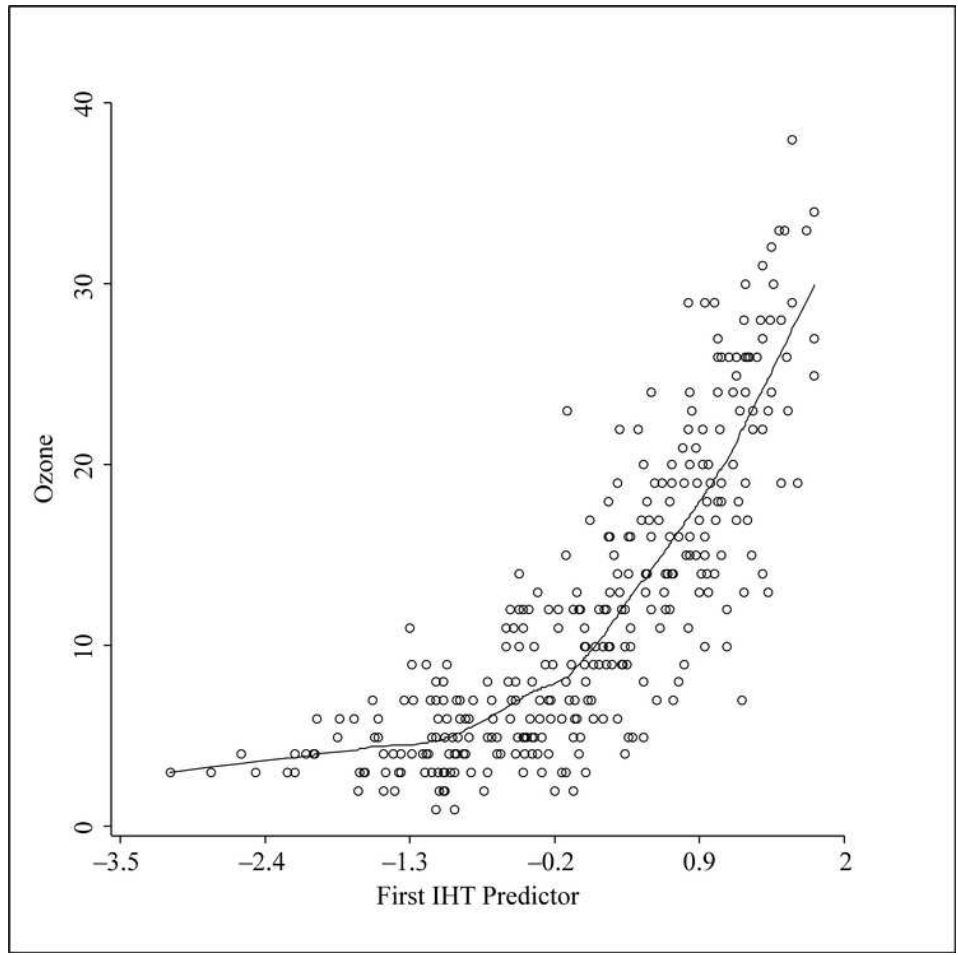

Fig. 1. *Scatterplot of ozone versus the first IHT predictor with a LOWESS smooth.*

Both tests are derived under the working coverage condition (C1), which is typically assumed in similar asymptotic developments in the literature. The working assumption is supported partly by the fact that IHT incorporates OLS and PHD in a way that does not evoke the constant variance condition (B), and is capable of discovering monotone and nonmonotone trends. Our experience indicates that IHT often works very well in picking up the patterns in a regression relation, so long as there is a nonzero vector to initiate the iteration process. Even if the coverage assumption does not hold, the general test still finds the significant vectors in the CMS, but the span of these vectors need not cover the CMS. In such cases, the tests should be viewed as a means of finding significant vectors in $\mathcal{S}_{\text{IHT}}$, which is a subspace of the CMS under the linear conditional mean condition (A). The issue of coverage is a fundamental and challenging one and deserves careful and in-depth investigation for IHT as well as other dimension reduction methods.

Cook and Critchley (2000) found that the CS automatically expands to incorporate regression outliers and mixtures. Consequently, they argued that the acknowledged sensitivity of CS methods like SIR and SAVE [see, e.g., Gather, Hilker and Becker (2002)] can be viewed as an advantage, since they have the ability to identify outliers and mixtures along with the main regression. In effect, methods for estimating the CS provide their own diagnostics. We conjecture that IHT is similarly self-diagnosing for outliers that affect the regression mean. Although we have not performed theoretical work to trace the diagnostic limits of IHT, various simulation results suggest that they might be fairly wide. For example, with four standard normal predictors, we generated 50 observations according to the linear model $Y = Z_1 + 0.2N(0,1)$, and then added a 51st observation with $Y = 6$ and corresponding $Z_j = 2$, $j = 1, \ldots, 4$. IHT estimated the dimension of the CMS to be 2, and the 3D summary plot clearly showed the linear mean structure and the outlier. Removal of the outlier resulted in a one-dimensional estimate of the CMS, as expected. Alternatively, we might deal with outliers by designing a robust version of IHT, replacing the sample moments by more robust estimators along the lines that Gather, Hilker and Becker (2001) used to investigate a robust version of SIR. This, however, is beyond the scope of the present paper.

The availability of these tests means that IHT is now a fully functioning methodology on a par with PHD. But, unlike PHD, it does not require the constant covariance condition (B) for either estimation or testing. In situations where PHD is applicable [conditions (C1)–(C3)], IHT automatically combines PHD with OLS, taking advantage of the ability of OLS to find linear trends in the mean function, and the ability of PHD to find nonlinear trends.

## APPENDIX: PROOFS OF LEMMAS

Throughout this section the identity matrix of dimension $p$ will be written as $I$ rather than $I_p$.

PROOF OF LEMMA 1. By definition,

$$\hat{\beta} = \widehat{\Sigma}^{-1/2} \hat{\sigma}^{-1} E_n (Z - \overline{Z})(Y - \overline{Y}). \tag{24}$$

Let us first expand $\widehat{\Sigma}^{-1/2}$. Note that

$$\begin{aligned}
\widehat{\Sigma} &= E_n(ZZ^T) - \overline{Z}\,\overline{Z}^T \\
&= E_n(ZZ^T) + O_p(n^{-1}) = I + E_n(ZZ^T - I) + O_p(n^{-1}),
\end{aligned}$$

where $E_n(ZZ^T - I)$ is of the order $O_p(n^{-1/2})$. We know that $\widehat{\Sigma}^{-1/2}$ must be of the form $I + A_n$ for some random matrix $A_n$ of the order $O_p(n^{-1/2})$. Therefore,

$$(I + A_n)^2 (I + E_n(ZZ^T - I)) = I.$$

The left-hand side is

$$I + E_n(ZZ^T - I) + 2A_n + O_p(n^{-1}).$$

Therefore $A_n = -E_n(ZZ^T - I)/2 + O_p(n^{-1})$ and

$$(25) \qquad \widehat{\Sigma}^{-1/2} = I - E_n(ZZ^T - I)/2 + O_p(n^{-1}).$$

By a similar argument one can show that

$$(26) \qquad \hat{\sigma}^{-1} = 1 - E_n(Y^2 - 1)/2 + O_p(n^{-1}).$$

It is easy to see that

$$(27) \qquad \begin{aligned} E_n(Z - \overline{Z})(Y - \overline{Y}) &= E_n(ZY) + O_p(n^{-1}) \\ &= \beta + E_n(ZY) - \beta + O_p(n^{-1}). \end{aligned}$$

Now substitute (25), (26) and (27) into (24) and expand the right-hand side of (24) to obtain expansion (4).

Next let us prove expansion (5). By definition,

$$(28) \qquad \widehat{H} = \widehat{\Sigma}^{-1/2} E_n[\hat{e}(Z - \overline{Z})(Z - \overline{Z})^T]\widehat{\Sigma}^{-1/2}.$$

We have already expanded $\widehat{\Sigma}^{-1/2}$. Now let us expand $E_n[\hat{e}(Z - \overline{Z})(Z - \overline{Z})^T]$. We have

$$E_n[\hat{e}(Z - \overline{Z})(Z - \overline{Z})^T] = E_n(\hat{e}ZZ^T) - \overline{Z}E_n(\hat{e}Z^T) - E_n(\hat{e}Z)\overline{Z}^T + O_p(n^{-1}).$$

Because $\overline{Z} = O_p(n^{-1/2})$, we need only expand $E_n(\hat{e}Z)$ so that the error is of the order $O_p(n^{-1/2})$. Note that

$$\begin{aligned} E_n(\hat{e}Z) &= E_n[(\hat{\sigma}^{-1}(Y - \overline{Y}) - \hat{\beta}^T\widehat{\Sigma}^{-1/2}(Z - \overline{Z}))Z] \\ &= \hat{\sigma}^{-1}E_n[(Y - \overline{Y})Z] - E_n[Z(Z - \overline{Z})^T]\widehat{\Sigma}^{-1/2}\hat{\beta}. \end{aligned}$$

It is easy to see that

$$\begin{aligned} \hat{\sigma}^{-1} &= 1 + O_p(n^{-1/2}), \\ E_n[(Y - \overline{Y})Z] &= \beta + O_p(n^{-1/2}), \\ E_n[Z(Z - \overline{Z})^T] &= I + O_p(n^{-1/2}), \\ \widehat{\Sigma}^{-1/2} &= I + O_p(n^{-1/2}), \\ \hat{\beta} &= \beta + O_p(n^{-1/2}). \end{aligned}$$

Therefore,
$$E_n(\hat{e}Z) = \beta - \beta + O_p(n^{-1/2}) = O_p(n^{-1/2}),$$

and consequently,

(29) $$E_n\hat{e}(Z - \overline{Z})(Z - \overline{Z})^T = E_n\hat{e}ZZ^T + O_p(n^{-1}).$$

We now expand the right-hand side so that the error is of the order $O_p(n^{-1})$. We have

(30) $$E_n\hat{e}ZZ^T = \hat{\sigma}^{-1}E_n[(Y - \overline{Y})ZZ^T] - E_n[\hat{\beta}^T\hat{\Sigma}^{-1/2}(Z - \overline{Z})(ZZ^T)].$$

The first term on the right-hand side is

(31) $$\begin{aligned}\hat{\sigma}^{-1}&E_n[(Y - \overline{Y})ZZ^T] \\ &= \hat{\sigma}^{-1}E_n[(Y - \overline{Y})(ZZ^T - I)] \\ &= \hat{\sigma}^{-1}E_n[Y(ZZ^T - I)] + O_p(n^{-1}) \\ &= (1 - E_n(Y^2 - 1)/2)E_n[Y(ZZ^T - I)] + O_p(n^{-1}).\end{aligned}$$

The second term on the right-hand side of (30) is expanded as
$$\begin{aligned}E_n[\hat{\beta}^T\hat{\Sigma}^{-1/2}&(Z - \overline{Z})(ZZ^T)] \\ &= E_n[\hat{\beta}^T\hat{\Sigma}^{-1/2}(Z - \overline{Z})(ZZ^T - I)] \\ &= E_n[\hat{\beta}^T\hat{\Sigma}^{-1/2}Z(ZZ^T - I)] + O_p(n^{-1}).\end{aligned}$$

The $(i,j)$th element of the $p \times p$ matrix on the right-hand side is
$$\sum_{k=1}^{p}(\hat{\Sigma}^{-1/2}\hat{\beta})_k E_n[Z_k(Z_iZ_j - \delta_{ij})],$$

where $(\hat{\Sigma}^{-1/2}\hat{\beta})_k$ is the $k$th element of the vector $\hat{\Sigma}^{-1/2}\hat{\beta}$ and $\delta_{ij}$ is the $(i,j)$th element of the $p$-dimensional identity matrix $I$. Because $Z$ has a standard multivariate normal distribution, the expectation of $Z_k(Z_iZ_j - \delta_{ij})$ is zero for any $i, j, k$. Therefore $E_n(Z_k(Z_iZ_j - \delta_{ij})) = O_p(n^{-1/2})$, and hence if we replace the $\hat{\Sigma}$ and $\hat{\beta}$ by $I$ and $\beta$, then the error incurred has the magnitude $O_p(n^{-1})$. It follows then that

(32) $$E_n[\hat{\beta}^T\hat{\Sigma}^{-1/2}(Z - \overline{Z})(ZZ^T)] = E_n[\beta^T Z(ZZ^T - I)] + O_p(n^{-1}).$$

Now substitute (31) and (32) into (30) to obtain
$$E_n(\hat{e}ZZ^T) = E_n[e(ZZ^T - I)] - \tfrac{1}{2}E_n(Y^2 - 1)E_n[Y(ZZ^T - I)] + O_p(n^{-1}).$$

However, note that
$$E[e(ZZ^T - I)] = E(eZZ^T) = H.$$

Hence,

$$E_n(\hat{e}ZZ^T) = H + E_n[e(ZZ^T - I) - H] - \tfrac{1}{2}E_n(Y^2 - 1)H + O_p(n^{-1}),$$

which, combined with (29), implies that

(33)
$$\begin{aligned}&E_n\hat{e}(Z - \overline{Z})(Z - \overline{Z})^T \\ &= H + E_n[e(ZZ^T - I) - H] - \tfrac{1}{2}E_n(Y^2 - 1)H + O_p(n^{-1}).\end{aligned}$$

Now substitute (25) and (33) into (28), and expand the right-hand side of (28) to obtain the desired expansion (5). □

PROOF OF LEMMA 2. (i) Since $\beta \notin \mathrm{Span}(H)$ and $H\beta, \ldots, H^{p-1}\beta$ belong to $\mathrm{Span}(H)$, $\beta \notin \mathrm{Span}\{H\beta, \ldots, H^{p-1}\beta\}$. Meanwhile, we know that

$$(\beta, H\beta, \ldots, H^{p-1}\beta)\Psi_0 = 0.$$

If the first row of $\Psi_0$ is not 0, then $\beta$ can be written as a linear combination of $H\beta, \ldots, H^{p-1}\beta$, which is a contradiction.

(ii) First, consider the case $\beta \perp \mathrm{Span}(H)$ (which includes the case $H = 0$). Then $B = (\beta, 0, \ldots, 0)$, $\mathrm{rank}(B) = 1$, and $\Psi_0$ is a $p$ by $p - 1$ matrix. Write

$$\Psi_0 = \begin{pmatrix} 0^T \\ \Phi_0 \end{pmatrix},$$

where $\Phi_0$ is a $p - 1$ by $p - 1$ matrix. Since $\Phi_0$ is an orthonormal matrix, its first row must contain a nonzero element.

Next, consider the case where $\beta$ is not orthogonal to $\mathrm{Span}(H)$. In this case $\mathrm{rank}(B) \geq 2$. Suppose first that $\mathrm{rank}(B) = 2$. Then $\Psi_0$ is a $p$ by $p - 2$ matrix. We claim that $H^2\beta \neq 0$. This is because if $H^2\beta = H(H\beta) = 0$, then $H\beta \perp \mathrm{Span}(H)$, but this implies $H\beta = 0$ since $H\beta$ belongs to $\mathrm{Span}(H)$. This means that $\beta \perp \mathrm{Span}(H)$, which is a contradiction. Hence

$$(H^2\beta, \ldots, H^{p-1}\beta) \neq 0.$$

Now suppose that the first row of $\Phi_0$ is 0 and write

$$\Phi_0 = \begin{pmatrix} 0 \\ \Lambda_0 \end{pmatrix},$$

where $\Lambda_0$ is a $p - 2$ by $p - 2$ matrix. Then $(H^2\beta, \ldots, H^{p-1}\beta)\Lambda_0 = 0$. In other words, the columns of $\Lambda_0$ are orthogonal to the rows of the matrix $(H^2\beta, \ldots, H^{p-1}\beta)$, which contains at least one nonzero row. Consequently the $p - 2$ columns of $\Lambda_0$ belong to a $(p - 3)$-dimensional space, so that they cannot be an orthogonal set. But this contradicts the fact that the columns of $\Psi_0$ are orthogonal.

Next, suppose that $\operatorname{rank}(B) = k > 2$. We first prove that $H\beta \in \operatorname{Span}(H^2\beta, \ldots, H^{p-1}\beta)$. From $\beta \notin \operatorname{Span}(H)$ it follows that the vectors $H^k\beta, \ldots, H^{p-1}\beta$ belong to the subspace spanned by the vectors $H\beta, \ldots, H^{k-1}\beta$, because otherwise we have, for some $j \in \{k, \ldots, p-1\}$ and some $c_1 \neq 0$,

$$H^j\beta = c_1\beta + c_2 H\beta + \cdots + c_k H^{k-1}\beta,$$

contradicting the assumption $\beta \notin \operatorname{Span}(H)$. By the same argument we can deduce that the vectors $H^k\beta, \ldots, H^{p-1}\beta$ must belong to the subspace spanned by $H^2\beta, \ldots, H^{k-1}\beta$. In particular,

$$H^k\beta = (H^2\beta, \ldots, H^{k-1}\beta)\delta$$

for some $\delta$ in $R^{k-2}$. Then

(34) $$H(H^{k-1}\beta - (H\beta, \ldots, H^{k-2}\beta)\delta) = 0.$$

In other words, the vector $H^{k-1}\beta - (H\beta, \ldots, H^{k-2}\beta)\delta$ is orthogonal to the rows, and hence columns, of $H$. However, both vectors in this difference belong to $\operatorname{Span}(H)$, and so we have

$$H^{k-1}\beta = (H\beta, \ldots, H^{k-2}\beta)\delta.$$

Consequently $H^{k-1}\beta$, and hence all the subsequent vectors $H^k, \ldots, H^{p-1}\beta$, belong to the space spanned by $H\beta, \ldots, H^{k-2}\beta$, which contradicts the assumption that $\operatorname{rank}(B) = k$.

However, if $H\beta$ belongs to $\operatorname{Span}(H^2\beta, \ldots, H^{p-1}\beta)$, then the matrix $(H^2\beta, \ldots, H^{p-1}\beta)$ has rank at least $k-1$, because we know that $H\beta, H^2\beta, \ldots, H^{k-1}\beta$ are linearly independent. Hence the solution space of the matrix $(H^2\beta, \ldots, H^{p-1}\beta)x = 0$ has dimension at most $(p-2) - (k-1) = p-k-1$. Now if the first two rows of $\Psi_0$ are zero, then there are $p-k$ orthogonal solutions to that equation, which is impossible. □

PROOF OF LEMMA 3. Since $\beta$ belongs to $\operatorname{Span}(H)$ it can be written as $\beta = H\eta$ for some $\eta$ in $\mathbb{R}^p$. First assume that $\operatorname{rank}(B) = 1$. We claim that $H^2\eta \neq 0$, otherwise $H\eta$ is orthogonal to $\operatorname{Span}(H)$, and must therefore be 0 because $H\eta$ belongs to $\operatorname{Span}(H)$. If the first row of $\Psi_0$ is 0, then

$$(H^2\eta, \ldots, H^p\eta)\Phi_0 = 0.$$

Therefore the $p-1$ columns of $\Phi_0$ are orthogonal to the rows of $(H^2\eta, \ldots, H^{p-1}\eta)$, which contains a nonzero row. But if so, the columns of $\Phi_0$ belong to a $(p-2)$-dimensional subspace of $\mathbb{R}^{p-1}$, and cannot be an orthogonal set— a contradiction.

Now suppose that $\operatorname{rank}(B) = k \geq 2$. We first prove that $\beta \in \operatorname{Span}(H\beta, \ldots, H^{p-1}\beta)$. Otherwise, by an argument similar to that used in Lemma 2, the vectors $H^k\beta, \ldots, H^{p-1}\beta$ all belong to the space spanned by $H\beta, \ldots, H^{k-1}\beta$. In particular, for some $\delta \in \mathbb{R}^{k-1}$,

$$H^k\beta = (H\beta, \ldots, H^{k-1}\beta)\delta,$$

which implies

$$H^{k+1}\eta = (H^2\eta, \ldots, H^k\eta)\delta.$$

But then, as we argued in the proof of Lemma 2, following display (34),

$$H^k\eta = (H\eta, \ldots, H^{k-1}\eta)\delta \quad \text{or} \quad H^{k-1}\beta = (\beta, \ldots, H^{k-2}\beta)\delta.$$

This implies that $H^{k-1}\beta$, and hence all its subsequent vectors $H^k\beta, \ldots, H^{p-1}\beta$, belong to the space spanned by $\beta, \ldots, H^{k-2}\beta$, contradicting the assumption that $\text{rank}(B) = k$.

That $\beta$ belongs to $\text{Span}(H\beta, \ldots, B^{p-1}\beta)$ implies that the matrix $(H\beta, \ldots, H^{p-1}\beta)$ has rank $k$, because $\beta, \ldots, H^{k-1}\beta$ are linearly independent. Therefore the equation $(H\beta, \ldots, H^{p-1}\beta)x = 0$ has at most $(p-1) - k = p - k - 1$ linearly independent solutions. However, if the first row of $\Psi_0$ is zero, then $(H\beta, \ldots, H^{p-1}\beta)x = 0$ has $p - k$ orthogonal solutions—a contradiction. $\square$

**Acknowledgments.** We are grateful to the Associate Editor and two referees for the helpful comments that led to a significant improvement of this paper.

SCHOOL OF STATISTICS  
UNIVERSITY OF MINNESOTA  
1994 BUFORD AVENUE  
ST. PAUL, MINNESOTA 55108  
USA  
E-MAIL: dennis@stat.umn.edu

DEPARTMENT OF STATISTICS  
PENNSYLVANIA STATE UNIVERSITY  
326 THOMAS BUILDING  
UNIVERSITY PARK, PENNSYLVANIA 16802  
USA  
E-MAIL: bing@stat.psu.edu